\documentclass[12pt,reqno]{amsart}
\usepackage{amsmath}
\usepackage{amssymb}
\usepackage{amsthm}

\newtheorem{theorem}{Theorem}[section]
\newtheorem{prop}[theorem]{Proposition}
\newtheorem{lem}[theorem]{Lemma}
\newtheorem*{theorema2}{Theorem A.2}
\newtheorem*{cor}{Corollary}
\theoremstyle{definition}
\newtheorem{defn}[theorem]{Definition}
\newtheorem*{defna1}{Definition A.1}
\theoremstyle{remark}
\newtheorem*{rem}{Remark}
\numberwithin{equation}{section}

\newcommand{\Gg}{\mathfrak{g}}    
\newcommand{\Gh}{\mathfrak{h}}
\newcommand{\Gq}{\mathfrak{q}}
\newcommand{\Gz}{\mathfrak{z}}

\begin{document}

\title[Quantum Heisenberg group algebra]
{${}^*$-representations of a quantum Heisenberg group algebra}

\author{Byung--Jay Kahng}
\date{}
\address{Department of Mathematics\\ University of Kansas\\
Lawrence, KS 66045}
\email{bjkahng@math.ukans.edu}

\begin{abstract}
In our earlier work, we constructed a specific non-compact quantum group
whose quantum group structures have been constructed on a certain twisted
group $C^*$-algebra.  In a sense, it may be considered as a ``quantum
Heisenberg group $C^*$-algebra''.  In this paper, we will find, up to
equivalence, all of its irreducible ${}^*$-representations.  We will point
out the Kirillov type correspondence between the irreducible representations
and the so-called ``dressing orbits''.  By taking advantage of its
comultiplication, we will then introduce and study the notion of {\em inner
tensor product representations\/}.  We will show that the representation
theory satisfies a ``quasitriangular'' type property, which does not appear
in ordinary group representation theory.
\end{abstract}
\maketitle

{\sc Introduction.} 
Recently in \cite{BJKp2}, we constructed a specific non-compact $C^*$-algebraic
quantum group, via deformation quantization of a certain non-linear Poisson--Lie
bracket on an exponential solvable Lie group.  The underlying $C^*$-algebra of
this quantum group has been realized as a twisted group $C^*$-algebra of a
nilpotent Lie group.  

 From its construction, it is reasonable to view it as a ``quantum Heisenberg
group $C^*$-algebra'' (This observation will be made a little clearer after the
first section.).  Focusing on its $C^*$-algebra structure, we study in this paper
its irreducible ${}^*$-representations.  It is not difficult to see that there
exists a Kirillov type, one-to-one correspondence between the irreducible
${}^*$-representations and the ``dressing orbits'' at the level of its
Poisson--Lie group counterpart.

 Since the object of our study is actually a Hopf $C^*$-algebra (i.\,e. quantum
group), we may use its comultiplication to define and study {\em inner tensor
product\/} of representations.  This is a generalization to quantum case of the
inner tensor product representations of an ordinary group.  We will show that
unlike in the case of an ordinary group and in the cases of many earlier examples
of non-compact quantum groups (e.\,g. \cite{VD}, \cite{SZ},\cite{Rf5}, \cite{Ld}),
the inner tensor product representations of our Hopf $C^*$-algebra satisfy a
certain ``quasitriangularity'' property.

 Here, we only study representation theory of the specific example of \cite{BJKp2}.
However, our earlier results (\cite{BJKp1,BJKp2}) imply that this quantum group
is just one example of a larger class of solvable quantum groups having twisted
group $C^*$-algebras or (more general) twisted crossed product algebras as
underlying $C^*$-algebras.  One of the main purposes of \cite{BJKp2} and this
paper is to present a case study, so that we are later able to develop a
procedure to construct and study more general class of locally compact quantum
groups.  Eventually, we wish to further develop a generalized orbit theory of
Kirillov type, which would then be used to study the harmonic analysis of the
locally compact quantum groups.  This will be our forthcoming project. 

 Since this paper is essentially a continuation of \cite{BJKp2}, we will
keep the same notation as in that paper.  Some of the notation is reviewed
in section 1.  This section is followed by Appendix, where we discuss the
classical counterparts (Poisson--Lie groups) to our quantum groups.  Main
purpose here is to calculate the dressing orbits on these Poisson--Lie groups.
Although the calculations are not difficult, these results could not be found
in the literatures.  Results here will be useful in our future study.

 In the second section, we discuss the representation theory of our example.
We find all the irreducible ${}^*$-representations up to equivalence.  We then
point out the Kirillov type one-to-one correspondence between these representations
and the dressing orbits.

 In the last section, we study inner tensor product representations.  By using
the quasitriangular quantum $R$-matrix operator obtained in \cite{BJKp2}, we
will show some interesting properties that do not appear in ordinary group
representation theory.

\section{The Hopf $C^*$-algebras}

 Our objects of study are the Hopf $C^*$-algebras $(A,\Delta)$ and $(\tilde{A},
\tilde{\Delta})$ constructed in \cite{BJKp2}.  As a $C^*$-algebra, $A$ is
isomorphic to a twisted group $C^*$-algebra.  That is, $A\cong C^*\bigl(H/Z,
C_{\infty}(\Gg/\Gq),\sigma\bigr)$, where $H$ is the $(2n+1)$ dimensional Heisenberg
Lie group (see Appendix) and $Z$ is the center of $H$.  Whereas, $\Gg=\Gh^*$ is
the dual space of the Lie algebra $\Gh$ of $H$ and $\Gq=\Gz^{\bot}$, for $\Gz
\subseteq\Gh$ corresponding to $Z$.  We denoted by $\sigma$ the twisting cocycle
for the group $H/Z$.  As constructed in \cite{BJKp2}, $\sigma$ is a continuous
field of cocycles $\Gg/\Gq\ni r\mapsto\sigma^r$, where
\begin{equation}\label{(sigma)}
\sigma^r\bigl((x,y),(x',y')\bigr)=\bar{e}\bigl[\eta_{\lambda}(r)\beta(x,y')\bigr].
\end{equation}
Following the notation of the previous paper, we used: $\bar{e}(t)=e^{(-2\pi i)t}$
and $\eta_{\lambda}(r)=\frac{e^{2\lambda r}-1}{2\lambda}$.  The elements $(x,y)$,
$(x',y')$ are group elements in $H/Z$.

 In \cite{BJKp2}, we showed that the $C^*$-algebra $A$ is a strict deformation
quantization (in the sense of Rieffel) of $C_{\infty}(\Gg)$, the commutative
algebra of continuous functions on $\Gg$ vanishing at infinity.  For convenience,
the deformation parameter $\hbar$ has been fixed ($\hbar=1$), which is the reason
why we do not see it in the definition of $A$.  When $\hbar=0$ (i.\,e. classical
limit), it turns out that $\sigma\equiv1$.  So $A_{\hbar=0}\cong C_{\infty}(G)$.
Throughout this paper (as in \cite{BJKp2}), we write $A=A_{\hbar=1}$.

 We could also construct a ``(regular) multiplicative unitary operator'' $U\in
{\mathcal B}({\mathcal H}\otimes{\mathcal H})$ associated with $A$, in the sense
of Baaj and Skandalis.  Thus $(A,\Delta)$ becomes a Hopf $C^*$-algebra, whose
comultiplication $\Delta$ is determined by $U$.  Similar realization as a twisted
group $C^*$-algebra also exists for the ``extended'' Hopf $C^*$--algebra $(\tilde{A},
\tilde{\Delta})$.  Again, its comultiplication $\tilde{\Delta}$ is determined by
a certain regular multiplicative unitary operator $\tilde{U}$.

 Actually, $(A,\Delta)$ is an example of a {\em locally compact quantum group\/},
equipped with the counit, antipode, and Haar weight, as constructed in \cite{BJKp2}.
We do not intend to give here the correct definition of a locally compact quantum
group, which is still at a primitive stage (But see \cite{KuV}, \cite{MN}, \cite{Wr7}
for some recent developments.).  Since the main goal of the present paper is in the
study of ${}^*$-representations of $A$, it would be rather sufficient to focus on
the Hopf $C^*$-algebra structure of $(A,\Delta)$.  For this reason, our preferred
terminology throughout this paper for $(A,\Delta)$ will be the ``Hopf $C^*$-algebra'',
although much stronger notion of the ``locally compact quantum group'' is still
valid.  Similar comments holds also for the extended Hopf $C^*$-algebra $(\tilde{A},
\tilde{\Delta})$.

 Before we continue our discussion, we insert the following Appendix.  It serves
three purposes: First, it briefly reviews the notation introduced in \cite{BJKp2}
and is being used here.  Second, it briefly gives a summary of the Poisson--Lie
group theory for those readers who are less familiar with the subject.  Third,
we calculate the dressing orbits on the Poisson--Lie groups which are classical
counterparts to our quantum groups.

\bigskip

\begin{center}
{\sc Appendix: Poisson--Lie groups, dressing actions and dressing orbits}
\end{center}

\noindent {\bf A.1. The Poisson--Lie groups.}
Let $H$ be the $(2n+1)$--dimensional Heisenberg Lie group such that the space
for the group is isomorphic to $\mathbb{R}^{2n+1}$ and the multiplication on it
is defined by
$$
 (x,y,z)(x',y',z')=\bigl(x+x',y+y',z+z'+\beta(x,y')\bigr),
$$
for $x,y,x',y'\in\mathbb{R}^n$ and $z,z'\in\mathbb{R}$.  Here $\beta(\ ,\ )$
is the usual inner product on $\mathbb{R}^n$, following the notation of 
\cite{BJKp2}.

 Consider also the extended Heisenberg Lie group $\tilde{H}$, with the group law
defined by
$$
(x,y,z,w)(x',y',z',w')=\bigl(x+e^{w}x',y+e^{-w}y',z+z'+(e^{-w})\beta(x,y'),
w+w'\bigr).
$$
The notation is similar as above, with $w,w'\in\mathbb{R}$.  This group contains
$H$ as a normal subgroup.

 In \cite{BJKp2}, we obtained the ``dual Poisson--Lie group'' $G$ of $H$.  It is
determined by the multiplication law:
$$
(p,q,r)(p',q',r')=(e^{\lambda r'}p+p',e^{\lambda r'}q+q',r+r'),
$$
while the dual Poisson--Lie group $\tilde{G}$ of $\tilde{H}$ is determined by the
multiplication law:
$$
(p,q,r,s)(p',q',r',s')=(e^{\lambda r'}p+p',e^{\lambda r'}q+q',r+r',s+s').
$$
Here $\lambda\in\mathbb{R}$ is a fixed constant, which determines a certain
non-linear Poisson structure on $G$ (or $\tilde{G}$) when $\lambda\ne0$.

 Although we do not explicitly mention the Poisson brackets here (see instead
\cite{BJKp2}), we can show indeed that $\tilde{H}$ and $\tilde{G}$ (similarly,
$H$ and $G$) are mutually dual Poisson--Lie groups.  For definition and some
important results on Poisson--Lie groups, see for example the article by Lu and
Weinstein \cite{LW} or the book by Chari and Pressley \cite{CP}.

 In \cite{BJKp2}, using the realization that the Poisson bracket on $G$ is a
non-linear Poisson bracket of the ``cocycle perturbation'' type, we have been
able to construct a quantum version of $G$:  The (non-commutative) Hopf
$C^*$--algebra $(A,\Delta)$, whose underlying $C^*$--algebra is a twisted group
$C^*$--algebra.  Similarly for $\tilde{G}$, we constructed the Hopf $C^*$--algebra
$(\tilde{A},\tilde{\Delta})$.  These are the main objects of study in \cite{BJKp2}
and in this paper.

 In this Appendix, we will study a special kind of a Lie group action of a
Poisson--Lie group on its dual Poisson--Lie group, called a ``dressing action''.
The orbits under the dressing action are the ``dressing orbits''.  We are
going to compute here the dressing orbits for our examples $\tilde{G}$ and
$\tilde{H}$ (and also $G$ and $H$).

\bigskip

\noindent {\bf A.2. Basic definitions: Dressing actions.}
Let $G$ be a Poisson--Lie group, let $G^*$ be its dual Poisson--Lie group,
and let $\Gg$ and $\Gg^*$ be the corresponding Lie algebras.  Together,
$(\Gg,\Gg^*)$ forms a Lie bialgebra.  On the vector space $\Gg\oplus\Gg^*$,
we can define a bracket operation by
\begin{align}
&\bigl[(X_1,\mu_1),(X_2,\mu_2)\bigr]  \notag\\
&=\bigl([X_1,X_2]_{\Gg}-\operatorname{ad}^*_{\mu_2}X_1
+\operatorname{ad}^*_{\mu_1}X_2,[\mu_1,\mu_2]_{\Gg^*}
+\operatorname{ad}^*_{X_1}\mu_2-\operatorname{ad}^*_{X_2}\mu_1\bigr),
\tag{A.1}
\end{align}
where $\operatorname{ad}^*_{X}\mu$ and $\operatorname{ad}^*_{\mu}X$ are,
respectively, the coadjoint representations of $\Gg$ on $\Gg^*$ and of
$\Gg^*$ on $\Gg=(\Gg^*)^*$.  This is a Lie bracket on $\Gg\oplus\Gg^*$,
which restricts to the given Lie brackets on $\Gg$ and $\Gg^*$.  We denote
the resulting Lie algebra by $\Gg\bowtie\Gg^*$, the {\em double Lie algebra\/}
\cite{LW}.

 Let $D=G\bowtie G^*$ be the connected, simply connected Lie group
corresponding to $\Gg\bowtie\Gg^*$.  There are homomorphisms of Lie groups
$$
G\hookrightarrow D\hookleftarrow G^*,
$$
lifting the inclusion maps of $\Gg$ and $\Gg^*$ into $\Gg\bowtie\Gg^*$.  Thus
we can define a product map $G\times G^*\to D$.  We will assume from now on
that the images of $G$ and $G^*$ are closed subgroups of $D$ and that the map
is a global diffeomorphism of $G\times G^*$ onto $D$.  In this case, we will
say that $D$ is a {\em double Lie group\/}.  In particular, each element of
$D$ has a unique expression $g\cdot\gamma$, for $g\in G$ and $\gamma\in G^*$.

 Suppose we are given a double Lie group $D=G\bowtie G^*$.  For $g\in G$ and
$\gamma\in G^*$, regarded naturally as elements in $D$, the product $\gamma\cdot
g$ would be factorized as
$$
\gamma\cdot g=g^{\gamma}\cdot\gamma^g,
$$
for some $g^{\gamma}\in G$ and $\gamma^g\in G^*$.  We can see without difficulty
that the map $\lambda:G^*\times G\to G$ defined by 
$$
\lambda_{\gamma}(g)=g^{\gamma}
$$
is a left action of $G^*$ on $G$.  Hence,
\begin{defna1}\label{dressing}
The map $\rho:G\times G^*\to G$ defined by
$$
\rho_{\gamma}(g)=g^{(\gamma^{-1})}
$$
is a right action, called the {\em dressing action\/} of $G^*$ on $G$.
\end{defna1}

\begin{rem}
Sometimes, the action $\lambda$ is called the left dressing action, while the
action $\rho$ is called the right dressing action.  It is customary to call the
right action $\rho$ the dressing action.  Semenov--Tian--Shansky \cite{Se} first
proved that the (right) dressing action of $G^*$ on $G$ is a Poisson action
(i.\,e. it preserves the respective Poisson structures).  The notion of dressing
action still exists (at least locally), even if the assumption that $G\times G^*
\cong D$ is not satisfied.  See \cite{LW}.
\end{rem}

 For any Lie algebra $\Gh$, it can always be regarded as a Lie bialgebra by
viewing its dual vector space $\Gg=\Gh^*$ as an abelian Lie algebra.  Then
the dressing action of $H$ on $G$ actually coincides with the coadjoint
action of $H$ on $\Gh^*\,(=\Gg\cong G)$.  In this sense, we may regard the
dressing action as a generalization of the coadjoint action.  This is the
starting point for the attempts to generalize the Kirillov's orbit theory,
and this point of view has been helpful throughout this paper.

 We conclude the subsection by stating the following important result by
Semenov--Tian--Shansky, which exhibits the close relationship between dressing
actions and the geometric aspects of Poisson--Lie groups.  This result appeared
in \cite{Se}, where he discusses the dressing actions in relation to the study
of complete integrable systems.  For the proof of the theorem, see \cite{Se}
or \cite{LW}.

\begin{theorema2}\label{a2}
The dressing action $G^*$ on $G$ is a Poisson action.  Moreover, the orbit of
the dressing action through $g\in G$ is exactly the symplectic leaf through
the point $g$ for the Poisson bracket on $G$.
\end{theorema2}

\bigskip

\noindent {\bf A.3. Dressing orbits for the Poisson--Lie groups $\tilde{G}$
and $\tilde{H}$.}
Let us consider our specific Poisson--Lie groups $\tilde{G}$ and $\tilde{H}$,
or equivalently, the Lie bialgebra $(\tilde{\Gg},\tilde{\Gh})$.  We will
construct here the double Lie group, dressing action, and dressing orbits.
Along the way, we will obtain the corresponding results for $G$ and $H$.

 By equation (A.1) and by using the Lie brackets on $\tilde{\Gg}$ and on
$\tilde{\Gh}$ (see \cite{BJKp2}), we can construct the double Lie algebra
$\tilde{\Gg}\bowtie\tilde{\Gh}$.  The space for it is $\tilde{\Gg}\oplus
\tilde{\Gh}$, on which the following Lie bracket is defined:
\small
\begin{align}
&\bigl[(p,q,r,s;x,y,z,w),(p',q',r',s';x',y',z',w')\bigr] \notag \\
&=\bigl(\lambda(r'p-rp')+(w'p-wp')+(r'y-ry'),\lambda(r'q-rq')+
(wq'-w'q)+(rx'-r'x),0,\notag \\
&\qquad(p'\cdot x-p\cdot x')+(q\cdot y'-q'\cdot y);(wx'-w'x)+
\lambda(rx'-r'x),(w'y-wy')+\lambda(ry'-r'y), \notag \\
&\qquad\beta(x,y')-\beta(x',y)+\lambda(p'\cdot x-p\cdot x')+
\lambda(q'\cdot y-q\cdot y'),0\bigr). \notag
\end{align}
\normalsize
We then calculate the corresponding Lie group $\tilde{D}$.  Using the notation
$\eta_{\lambda}(r)=\frac{e^{2\lambda r}-1}{2\lambda}$ (the function introduced
in Definition 2.3 of \cite{BJKp2}), the group $\tilde{D}$ is given by the
following multiplication law:
\small
\begin{align}
&(p,q,r,s;x,y,z,w)(p',q',r',s';x',y',z',w')  \notag \\
&=\bigl(e^{\lambda r'}p+e^{-w}p'+e^{-\lambda r'}\eta_{\lambda}(r')y,
e^{\lambda r'}q+e^{w}q'-e^{-\lambda r'}\eta_{\lambda}(r')x,r+r',
\notag \\
&\qquad s+s'+(e^{-\lambda r'-w})p'\cdot x-(e^{-\lambda r'+w})q'\cdot y-
\eta_{\lambda}(-r')\beta(x,y);e^{-\lambda r'}x+e^{w}x',e^{-\lambda r'}y
+e^{-w}y',  \notag \\
&\qquad z+z'+(e^{-\lambda r'-w})\beta(x,y')+\lambda(e^{-\lambda r'-w})p'
\cdot x+\lambda(e^{-\lambda r'+w})q'\cdot y+\lambda\eta_{\lambda}(-r')
\beta(x,y),w+w'\bigr). \notag
\end{align}
\normalsize

 If we identify $(p,q,r,s)\in\tilde{G}$ with $(p,q,r,s;0,0,0,0)$ and
$(x,y,z,w)\in{\tilde{H}}$ with $(0,0,0,0;x,y,z,w)$, it is easy to see that
$\tilde{G}$ and $\tilde{H}$ are (closed) Lie subgroups of $\tilde{D}$. 
This defines a global diffeomorphism of $\tilde{G}\times\tilde{H}$ onto
$\tilde{D}$, since any element $(p,q,r,s;x,y,z,w)$ of $\tilde{D}$ can be
written as
$$
(p,q,r,s;x,y,z,w)=(p,q,r,s;0,0,0,0)(0,0,0,0;x,y,z,w).
$$
In other words, $\tilde{D}$ is the {\em double Lie group\/} of $\tilde{G}$
and $\tilde{H}$, which will be denoted by $\tilde{G}\bowtie\tilde{H}$.
Meanwhile, if we consider only the $(p,q,r)$ and the $(x,y,z)$ variables,
we obtain in the same way the double Lie group $D=G\bowtie H$ of $G$ and $H$.

 By Definition A.1, the {\em dressing action\/} of $\tilde{H}$ on $\tilde{G}$
is:
\begin{align}
\bigl(\rho(x,y,z,w)\bigr)(p,q,r,s)&=\bigl(e^{w}p-e^{-\lambda r+w}\eta
_{\lambda}(r)y,e^{-w}q+e^{-\lambda r-w}\eta_{\lambda}(r)x,r,
\notag \\
&\qquad s-e^{-\lambda r}p\cdot x+e^{-\lambda r}q\cdot y
+e^{-2\lambda r}\eta_{\lambda}(r)\beta(x,y)\bigr).  \notag
\end{align}
The dressing orbits, which by Theorem A.2 are the symplectic leaves
in $\tilde{G}$ for its (non-linear) Poisson bracket, are:
\begin{itemize}
\item $\tilde{\mathcal O}_s=\{(0,0,0,s)\}$, when $(p,q,r)=(0,0,0)$.
\item $\tilde{\mathcal O}_{p,q}=\{(ap,\frac1a q,0,c):a>0,c\in\mathbb{R}\}$,
when $r=0$ but $(p,q)\ne(0,0)$.
\item $\tilde{\mathcal O}_{r,s}=\bigl\{\bigl(a,b,r,s-\frac1{\eta_{\lambda}
(r)}a\cdot b\bigr):(a,b)\in\mathbb{R}^{2n}\bigr\}$, when $r\ne0$.
\end{itemize}
Here $a\cdot b$ denotes the inner product.  The $\tilde{\mathcal O}_s$ are
1--point orbits, the $\tilde{\mathcal O}_{p,q}$ are 2--dimensional orbits,
and the $\tilde{\mathcal O}_{r,s}$ are $2n$--dimensional orbits.

 Similarly, if we only consider the $(p,q,r)$ and the $(x,y,z)$ variables,
we obtain the expression for the dressing action of $H$ on $G$ as follows:
$$
\bigl(\rho(x,y,z)\bigr)(p,q,r)=\bigl(p-e^{-\lambda r}\eta_{\lambda}(r)y,
q+e^{-\lambda r}\eta_{\lambda}(r)x,r\bigr).
$$
So the dressing orbits in $G$ are:
\begin{itemize}
\item (1--point orbits): ${\mathcal O}_{p,q}=\{(p,q,0)\}$, when 
$r=0$.
\item ($2n$--dimensional orbits): ${\mathcal O}_r=\{(a,b,r):(a,b)
\in\mathbb{R}^{2n}\}$, when $r\ne0$.
\end{itemize}

 Meanwhile, to calculate the dressing action of $\tilde{G}$ on $\tilde{H}$,
it is convenient to regard $\tilde{D}$ as the double Lie group $\tilde{H}
\bowtie\tilde{G}$ of $\tilde{H}$ and $\tilde{G}$.  Indeed, there exists a
global diffeomorphism between $\tilde{H}\times\tilde{G}$ and $\tilde{D}$
defined by
\small
\begin{align}
&\bigl((x,y,z,w),(p,q,r,s)\bigr)\mapsto(0,0,0,0;x,y,z,w)(p,q,r,s,0,0,0,0)
\notag \\
&=\bigl(e^{-w}p+e^{-\lambda r}\eta_{\lambda}(r)y,e^{w}q-e^{-\lambda r}
\eta_{\lambda}(r)x,r,s+(e^{-\lambda r-w})p\cdot x-(e^{-\lambda r+w})q\cdot y
+\eta_{\lambda}(r)\beta(x,y);  \notag \\
&\qquad e^{-\lambda r}x,e^{-\lambda r}y,z+\lambda(e^{-\lambda r-w})p\cdot x
+\lambda(e^{-\lambda r+w})q\cdot y-\lambda\eta_{\lambda}(r)\beta(x,y),w\bigr). 
\notag
\end{align}
\normalsize
Similarly, we can show that $D$ is the double Lie group of $H$ and $G$.

 Using this characterization of the double Lie group $\tilde{H}\bowtie
\tilde{G}$, the dressing action of $\tilde{G}$ on $\tilde{H}$ is obtained
by Definition A.1.  That is,
\begin{align}
&\bigl(\rho(p,q,r,s)\bigr)(x,y,z,w)  \notag \\
&=\bigl(e^{-\lambda r}x,e^{-\lambda r}y,z+\lambda(e^{-\lambda r})
p\cdot x +\lambda(e^{-\lambda r})q\cdot y-\lambda(e^{-2\lambda r})
\eta_{\lambda}(r)\beta(x,y),w\bigr). \notag
\end{align}
The dressing orbits in $\tilde{H}$ are:
\begin{itemize}
\item (1--point orbits): $\tilde{\mathcal O}_{z,w}=\{(0,0,z,w)\}$, when 
$(x,y)=(0,0)$.
\item (2--dimensional orbits): $\tilde{\mathcal O}_{x,y,w}=\{(\alpha x,
\alpha y,\gamma,w):\alpha>0,\gamma\in\mathbb{R}\}$, when $(x,y)\ne(0,0)$.
\end{itemize}

 Similarly, the dressing action of $G$ on $H$ is given by
$$
\bigl(\rho(p,q,r)\bigr)(x,y,z)=\bigl(e^{-\lambda r}x,e^{-\lambda r}y,
z+\lambda(e^{-\lambda r})p\cdot x-\lambda(e^{-\lambda r})q\cdot y
-\lambda(e^{-2\lambda r})\eta_{\lambda}(r)\beta(x,y)\bigr),
$$
and the dressing orbits in $H$ are:
\begin{itemize}
\item (1--point orbits): ${\mathcal O}_{z}=\{(0,0,z)\}$, when $(x,y)
=(0,0)$.
\item (2--dimensional orbits): ${\mathcal O}_{x,y}=\{(\alpha x,\alpha y,
\gamma):\alpha>0,\gamma\in\mathbb{R}\}$, when $(x,y)\ne(0,0)$.
\end{itemize}

\begin{rem}
As pointed out earlier, the dressing action is usually regarded as a
generalization of the coadjoint action.  In the present case, we can
see easily that the dressing orbits in $\tilde{H}$ are exactly the
coadjoint orbits in $\tilde{\Gh}\cong\tilde{H}$.  This illustrates
the point that the Poisson bracket on $\tilde{H}$ is just the linear
Poisson bracket.  On the other hand, for $\tilde{G}$, which has a
non-linear Poisson bracket, this is no longer the case.  The orbits
$\tilde{\mathcal O}_{r,s}$ are different from the coadjoint orbits
in $\tilde{\Gg}$.  Nevertheless, we can still see close resemblance.
\end{rem}

\section{Quantum Heisenberg group algebra representations}

 From its construction, we can see that $(A,\Delta)$ is at the same time a
``quantum $C_{\infty}(G)$'' and ``quantum $C^*(H)$''.  In the previous paper
\cite{BJKp2}, the first viewpoint has been exploited:  As we already mentioned,
$(A,\Delta)$ has been constructed as a deformation quantization of $C_{\infty}(G)$.
The construction of the counit, antipode, and Haar weight for $(A,\Delta)$ all
comes from the corresponding structures on $G$.

 In this article, we wish to focus our attention to the second viewpoint.  To
see this, recall that $A\cong C^*\bigl(H/Z,C_{\infty}(\Gg/\Gq),\sigma\bigr)$,
where the twisting cocycle $\sigma$ for $H/Z$ is defined as in \eqref{(sigma)}.
Since we may put $\eta_{\lambda}(r)=r$ for $\lambda=0$, it is a simple exercise
using Fourier inversion theorem that $A\cong C^*(H)$ when $\lambda=0$ (The reader
may refer to the article \cite{BJKp1} or \cite{Rf2} for the definition of a
twisted group $C^*$-algebra.).  For this reason, we will on occasion call $(A,
\Delta)$ the {\em quantum Heisenberg group $C^*$-algebra\/}.

\begin{rem}
The notion of the ``quantum Heisenberg group ($C^*$-)algebra'' introduced above
is different from the notion of the ``quantum Heisenberg algebra'' used in some
physics literatures \cite{GF}, \cite{Ke}.  They are different as algebras.
Another significant distinction is that ours is equipped with a (non-cocommutative)
comultiplication, while the other one does not consider any coalgebra structure.
Hence the slight difference in the choice of the terminologies.
\end{rem}

 The study of the ${}^*$-representations of $(A,\Delta)$ (and also of $(\tilde{A},
\tilde{\Delta})$) will be a generalization of the study of the Heisenberg group
representation theory (which, by a standard result, is equivalent to the
${}^*$-representation theory of $C^*(H)$).  The topic itself is of interest to us
(providing us with the properties like ``quasitriangularity'').  But the success
of the (ordinary) Heisenberg group representation theory in various applications
also suggests that this is a worthwhile topic to develop.

 Before we begin our discussion, let us fix our terminology.  Suppose $(A,\Delta)$
is a general Hopf $C^*$-algebra (in the sense of \cite{Va}, \cite{BS}).  By a
{\em representation\/} of $(A,\Delta)$, we will just mean a non-degenerate
${}^*$-representation of the $C^*$-algebra $A$ on a certain Hilbert space.

 On the other hand, by a {\em coaction\/} of $(A,\Delta)$ on a $C^*$-algebra $B$,
we will mean a non-degenerate ${}^*$-homomorphism $\delta_B:B\to M(B,A)$ such that 
$$
(\operatorname{id}_B\otimes\Delta)\delta_B=(\delta_B\otimes\operatorname{id}_A)
\delta_B.
$$
Here $M(B,A)$ is the set $\{x\in M(B\otimes A):x(1_{M(B)}\otimes A)+(1_{M(B)}
\otimes A)x\subseteq B\otimes A\}$, which is a $C^*$-subalgebra of the multiplier
algebra $M(B\otimes A)$.  Similarly, a {\em (unitary) corepresentation\/} of the
Hopf $C^*$-algebra $(A,\Delta)$ on a Hilbert space ${\mathcal H}$ is a unitary
$\Pi\in M\bigl({\mathcal K}({\mathcal H})\otimes A\bigr)$ such that:
\begin{equation}\label{(corep)}
(\operatorname{id}\otimes\Delta)(\Pi)=\Pi_{12}\Pi_{13}.
\end{equation}
Here $\Pi_{12}$ is understood as an element in $M\bigl({\mathcal K}({\mathcal H})
\otimes A\otimes A\bigr)$ such that it acts as $\Pi$ on the first and the second
variables while the remaining variable is unchanged.  The notation $\Pi_{13}$
is understood in the similar manner.

\begin{rem}
Corepresentations of $(A,\Delta)$ are actually the ``representations of the
coalgebra structure on $(A,\Delta)$''.  So in many articles on quantum groups,
they are often called ``(unitary) representations of the locally compact quantum
group $(A,\Delta)$''.  In particular, representation theory in this sense of
{\em compact quantum groups\/} \cite{Wr2}, which are themselves Hopf $C^*$-algebras,
have been neatly studied by Woronowicz in \cite{Wr3}.  However, note that we will
use the terminologies ``representations'' and ``corepresentations'' of a Hopf
$C^*$-algebra in the sense defined above.  This would make things a little
simpler.  Moreover, this is closer to the spirit of this paper, trying to view
our $(A,\Delta)$ as a quantum Heisenberg group $C^*$-algebra.
\end{rem} 

 Recall that in our case, the $C^*$-algebra $A$ is isomorphic to the twisted
group $C^*$-algebra $C^*\bigl(H/Z,C_{\infty}(\Gg/\Gq),\sigma\bigr)$.  So by
slightly modifying Theorem 3.3 and Proposition 3.4 of \cite{BuS}, we are able
to obtain the representations of $A$ from the so-called ``representing pairs''
$(\mu,Q^{\sigma})$.  Such a pair $(\mu,Q^{\sigma})$ consists of a nondegenerate
representation $\mu$ of $C_{\infty}(\Gg/\Gq)$ and a generalized projective
representation $Q^{\sigma}$ of $H/Z$, satisfying the following property:
\begin{equation}\label{(genproj)}
Q^{\sigma}_{(x,y)}Q^{\sigma}_{(x',y')}=\mu\bigl(\sigma\bigl((x,y),(x',y')
\bigr)\bigr)Q^{\sigma}_{(x+x',y+y')},
\end{equation}
for $(x,y),(x',y')\in H/Z$.

 Given a representing pair $(\mu,Q^{\sigma})$, we can construct its ``integrated
form''.  On the dense subspace of Schwartz functions, it reads:
\begin{equation}\label{(integratedform)}
\pi(f)=\int_{H/Z}\mu\bigl(f(x,y;\cdot)\bigr)Q^{\sigma}_{(x,y)}\,dxdy.
\end{equation}
By natural extension to the $C^*$-algebra level (Due to the amenability of the
group involved, there is no ambiguity.  See \cite{BJKp1}.), we obtain in this way
a representation of $A$.

\begin{rem}
Let us from now on denote by ${\mathcal A}$ the dense subspace $S_{3c}(H/Z
\times\Gg/\Gq)$ of $A$, which is the space of Schwartz functions in the
$(x,y;r)$ variables having compact support in the $r\,(\in\Gg/\Gq)$ variable.
This is a dense subalgebra (under the twisted convolution) of our twisted group
$C^*$-algebra $A$, and it has been used throughout \cite{BJKp2} (However,
we should point out that our usage of ${\mathcal A}$ is slightly different
from that of \cite{BJKp2}: There, ${\mathcal A}$ is contained in $S(\Gg)$, while
at present we view it as functions in the $(x,y;r)$ variables.   Nevertheless,
they can be regarded as the same if we consider these functions as operators
contained in our $C^*$-algebra.).  Similarly for $\tilde{A}$, we will consider
the dense subalgebra $\tilde{\mathcal A}$ of Schwartz functions in the $(x,y,r,w)$
variables having compact support in the $r$ and $w$.
\end{rem}

 To find irreducible representations of $A$, let us look for some representing
pairs $(\mu,Q^{\sigma})$ consisting of irreducible $\mu$ and $Q^{\sigma}$.
Irreducible representations of the commutative algebra $C_{\infty}(\Gg/\Gq)$
are just the pointwise evaluations at $r\in\Gg/\Gq$.  So let us fix $r\in
\Gg/\Gq$ and the corresponding 1-dimensional representation $\mu$ of $C_{\infty}
(\Gg/\Gq)$, given by $\mu(v)=v(r)$, for $v\in C_{\infty}(\Gg/\Gq)$.  Then the 
condition for $Q^{\sigma}$ becomes:
\begin{align}
Q^{\sigma}_{(x,y)}Q^{\sigma}_{(x',y')}&=\sigma\bigl((x,y),(x',y');r\bigr)
Q^{\sigma}_{(x+x',y+y')}  \notag \\
&=\sigma^r\bigl((x,y),(x',y')\bigr)Q^{\sigma}_{(x+x',y+y')}. \label{(projrepn)}
\end{align}
That is, $Q^{\sigma}$ satisfies the condition for an (ordinary) projective
representation of $H/Z$ with respect to the ordinary $\mathbb{T}$--valued
cocycle $\sigma^r$.

 Using $\sigma^r$, we may define an extension group $E$ of $H/Z$.  Its
underlying space is $H/Z\times\mathbb{T}$ and its multiplication law is
given by
\begin{align}
(x,y;{\theta})(x',y';{\theta}')&=\bigl(x+x',y+y';{\theta}{\theta}'\sigma^r
\bigl((x,y),(x',y')\bigr)\bigr)  \notag\\
&=\bigl(x+x',y+y';{\theta}{\theta}'\bar{e}\bigl[\eta_{\lambda}(r)\beta(x,y')
\bigl]\bigr).   \notag
\end{align}
Standard theory tells us that the unitary projective representations of $(H/Z,
\sigma^r)$ come from the unitary group representations of $E$, which is easier
to study:  The next lemma gives us all the irreducible unitary representations
of $E$, up to equivalence.

\begin{lem}
Let $E$ be the extension group of $(H/Z,\sigma^r)$ as defined above.  Then its
irreducible unitary representations are equivalent to one of the following:
\begin{itemize}
\item For each $(p,q)\in\mathbb{R}^{2n}$, there exists a 1--dimensional
representation $Q_{p,q}$ defined by
$$
Q_{p,q}(x,y;\theta)=\bar{e}(p\cdot x+q\cdot y).
$$
\item There also exists a (unique) infinite dimensional representation $Q_r$
on $L^2(\mathbb{R}^n)$ defined by
$$
\bigl(Q_r(x,y;\theta)\xi\bigr)(u)=\theta\bar{e}\bigl[\eta_{\lambda}(r)
\beta(u,y)\bigr]\xi(u+x).
$$
\end{itemize}
\end{lem}

\begin{proof}
Observe that $E$ is a semi-direct product of two abelian groups $X=\{(x,0,1):
x\in\mathbb{R}^n\}$ and $Y\times\mathbb{T}=\{(0,y,\theta):y\in\mathbb{R}^n,
\theta\in\mathbb{T}\}$.  So by using Mackey analysis, every irreducible
representation of $E$ is obtained as an ``induced representation'' \cite{FD},
\cite{Tay}, \cite{Rfm}.

Since $[E,E]=\mathbb{T}$, we have: $E/[E,E]=X\times Y$, which is abelian.  So
all the irreducible representations of $E/[E,E]$ are one-dimensional.  By
lifting from these 1-dimensional representations, we obtain the (irreducible)
representations $\{Q_{p,q}\}_{(p,q)\in\mathbb{R}^{2n}}$ of $E$ that are trivial
on the commutator $[E,E]$.

The infinite dimensional representation $Q_r$ is the induced representation
$\operatorname{Ind}_{Y\times\mathbb{T},\chi}^E$, where $\chi$ is the
representation of $Y\times\mathbb{T}$ defined by: $\chi(y,\theta)=\theta$.
It turns out (by using standard Mackey theory) that $\{Q_{p,q}\}_{(p,q)\in
\mathbb{R}^{2n}}$ and $Q_r$ exhaust all the irreducible representations of
$E$, up to equivalence.
\end{proof}

 We are now able to find the irreducible projective representations of $(H/Z,
\sigma)$.  Check equation \eqref{(projrepn)} and we obtain the following
representing pairs consisting of irreducible $\mu$ and $Q^{\sigma}$.
\begin{enumerate}
\item (When $r=0\in\Gg/\Gq$):  For each $(p,q)\in\mathbb{R}^{2n}$, there is a
pair $(\mu,Q^{\sigma})$ given by
\begin{itemize}
\item $\mu(v)=v(0)$, $v\in C_{\infty}(\Gg/\Gq)$.
\item $Q^{\sigma}(x,y)=\bar{e}(p\cdot x+q\cdot y),\quad (x,y)\in H/Z$.
\end{itemize}
\item (When $r\ne0\in\Gg/\Gq$):  There is a pair $(\mu,Q^{\sigma})$ given by
\begin{itemize}
\item $\mu(v)=v(r)$, $v\in C_{\infty}(\Gg/\Gq)$.
\item On $L^2(\mathbb{R}^n)$, 
$$
\bigl(Q^{\sigma}(x,y)\xi\bigr)(u)=\bar{e}\bigl[\eta_{\lambda}(r)\beta(u,y)\bigr]
\xi(u+x),\quad (x,y)\in H/Z.
$$
\end{itemize}
\end{enumerate}

 Therefore, we obtain the following proposition.  Observe the similarity
between this result and the representation theory of the Heisenberg group
$H$ or the Heisenberg group $C^*$-algebra $C^*(H)$ (This is not surprising,
given the remark at the beginning of this section.).

\begin{prop}\label{repA}
Consider the twisted convolution algebra ${\mathcal A}$.  Its irreducible
representations are equivalent to one of the following representations,
which have been obtained by integrating the representing pairs $(\mu,
Q^{\sigma})$ of the preceding paragraph.
\begin{itemize}
\item For $(p,q)\in\mathbb{R}^{2n}$, there is a 1-dimensional 
representation ${\pi}_{p,q}$ of ${\mathcal A}$, defined by
$$
\pi_{p,q}(f)=\int f(x,y,0)\bar{e}(p\cdot x+q\cdot y)\,dxdy.
$$ 
\item For $r\in\mathbb{R}$, there is a representation ${\pi}_r$ of
${\mathcal A}$, acting on the Hilbert space ${\mathcal H}_r=L^2(\mathbb{R}^n)$
and is defined by
$$
\bigl(\pi_r(f)\xi\bigr)(u)=\int f(x,y,r)\bar{e}\bigl[\eta_{\lambda}(r)\beta(u,y)
\bigr]\xi(u+x)\,dxdy.
$$
\end{itemize}
Since ${\mathcal A}$ is a dense subalgebra of our $C^*$-algebra $A$, we thus
obtain all the irreducible representations (up to equivalence) of $A$ by
naturally extending these representations.  We will use the same notation,
$\pi_{p,q}$ and $\pi_r$, for the representations of $A$ constructed in this way.
\end{prop}

 Let us now consider the representations of the $C^*$-algebra $\tilde{A}$.
They are again obtained from representations of the dense subalgebra $\tilde
{\mathcal A}$, which have been identified with the twisted convolution algebra
of functions in the $(x,y,r,w)$ variables (where $(x,y,w)\in\tilde{H}/Z$ and
$r\in\tilde{\Gg}/{\Gz^{\bot}}$) having compact support in the $r$ and $w$
variables.  We may employ the same argument as above to find (up to equivalence)
the irreducible representations of $\tilde{\mathcal A}$.  The result is given
in the following proposition:

\begin{prop}\label{repAtilde}
The irreducible representations of $\tilde{A}$ are obtained by naturally
extending the following irreducible representations of the dense subalgebra
$\tilde{\mathcal A}$.
\begin{itemize}
\item For $s\in\mathbb{R}$, there is a 1-dimensional representation
$\tilde{\pi}_s$ defined by
$$
\tilde{\pi}_{s}(f)=\int f(x,y,0,w)\bar{e}(sw)\,dxdydw.
$$
\item For $(p,q)\in\mathbb{R}^{2n}$, there is a representation $\tilde
{\pi}_{p,q}$ acting on the Hilbert space $\tilde{\mathcal H}_{p,q}=L^2
(\mathbb{R})$ defined by
$$
\bigl(\tilde{\pi}_{p,q}(f)\zeta\bigr)(d)=\int f(x,y,0,w)\bar{e}(e^d p\cdot x+e^{-d} 
q\cdot y)\zeta(d+w)\,dxdydw.
$$
\item For $(r,s)\in\mathbb{R}^2$, there is a representation $\tilde{\pi}
_{r,s}$ acting on the Hilbert space $\tilde{\mathcal H}_{r,s}=L^2
(\mathbb{R}^n)$ defined by
$$
\bigl(\tilde{\pi}_{r,s}(f)\xi\bigr)(u)=\int f(x,y,r,w)\bar{e}(sw)\bar{e}\bigl
[\eta_{\lambda}(r)\beta(u,y)\bigr](e^{-\frac{w}2})^n\xi(e^{-w}u+e^{-w}x)
\,dxdydw.
$$
\end{itemize}
We will use the same notation, $\tilde{\pi}_s$, $\tilde{\pi}_{p,q}$ and
$\tilde{\pi}_{r,s}$, for the corresponding representations of $\tilde{A}$.
\end{prop}

\begin{proof}
As before, let us first fix $r\in{\Gg}/{\Gz^{\bot}}$.  Look for the
irreducible projective representations of $\tilde{H}/Z$, with respect to
the ($\mathbb{T}$-valued) cocycle for $\tilde{H}/Z$ defined by
$$
\tilde{\sigma}^r:\bigl((x,y,w),(x',y',w')\bigr)\mapsto\bar{e}
\bigl[e^{-w}\eta_{\lambda}(r)\beta(x,y')\bigr].
$$
To do this, we consider the extension group $\tilde{E}$ of $\tilde{H}/Z$,
whose underlying space is $\tilde{H}/Z\times\mathbb{T}$ and whose
multiplication law is given by
$$
(x,y,w;\theta)(x',y',w';\theta')=\bigl(x+e^wx',y+e^{-w}y',w+w';
\theta\theta'\bar{e}\bigl[e^{-w}\eta_{\lambda}(r)\beta(x,y')\bigr]\bigr).
$$

Again, all the irreducible representations of $\tilde{E}$ are obtained by
``inducing''.  Up to equivalence, they are:
\begin{itemize}
\item For each $s\in\mathbb{R}$, there exists a 1--dimensional representation
$\tilde{Q}_s$ defined by
$$
\tilde{Q}_s(x,y,w;\theta)=\bar{e}(sw).
$$
\item For each $(p,q)\in\mathbb{R}^{2n}$, there exists a representation
$\tilde{Q}_{p,q}$ on $L^2(\mathbb{R})$ defined by
$$
\bigl(\tilde{Q}_{p,q}(x,y,w;\theta)\zeta\bigr)(d)=\bar{e}(e^dp\cdot x+e^{-d}q\cdot y)
\zeta(d+w).
$$
\item For $s\in\mathbb{R}$, there exists an infinite dimensional irreducible
representation $\tilde{Q}_{r,s}$ on $L^2(\mathbb{R}^n)$ defined by
$$
\bigl(\tilde{Q}_{r,s}(x,y,w;\theta)\xi\bigr)(u)=\theta\bar{e}(sw)\bar{e}\bigl[
\eta_{\lambda}(r)\beta(u,y)\bigr](e^{-\frac{w}2})^n\,\xi(e^{-w}u+e^{-w}x).
$$
\end{itemize}

 Vary $r\in{\Gg}/{\Gz^{\bot}}$ and check the compatibility condition just
like \eqref{(projrepn)}, to find the appropriate representing pairs.  Then
the integrated form of these pairs will give us the irreducible representations
of $\tilde{\mathcal A}$, which are stated in the proposition.
\end{proof}

 As in the case of Proposition \ref{repA} (for the $C^*$-algebra $A$),
we can see clearly the similarity between the result of Proposition
\ref{repAtilde} and the representation theory of the ordinary group
$C^*$-algebra $C^*(\tilde{H})$.  Indeed, except when we study later the
notion of ``inner tensor product representations'' (taking advantage of the
Hopf structures of $A$ and $\tilde{A}$), the representation theories of $A$
and $\tilde{A}$ are very similar to those of $C^*(H)$ and $C^*(\tilde{H})$.

 In this light, it is interesting to observe that the irreducible
representations of $A$ and $\tilde{A}$ are in one-to-one correspondence
with the dressing orbits (calculated in Appendix) in $G$ and $\tilde{G}$,
respectively.  To emphasize the correspondence, we used the same subscripts
for the orbits and the related irreducible representations.  In this paper,
we will point out this correspondence only.  However, it is still true
(as in the ordinary Lie group representation theory) that orbit analysis
sheds some helpful insight into the study of quantum group representations.

 In the following, we give a useful result about the irreducible
representations of $\tilde{A}$ and those of $A$.  This has been 
motivated by the orbit analysis, and it is an analog of a similar
result for the group representations of $\tilde{H}$ and $H$.  See
the remark following the proposition.

\begin{defn}\label{restriction}
Suppose we are given a representation $\tilde{\pi}$ of $\tilde{A}$.  Since
it is essentially obtained from a representation $\tilde{Q}$ of $\tilde{E}$,
we may consider its restriction $\tilde{Q}|_E$ to $E$.  Let us denote by
$\tilde{\pi}|_A$ the representation of $A$ corresponding to the representation
$\tilde{Q}|_E$ of $E$.  In this sense, we will call $\tilde{\pi}|_A$ the
{\em restriction\/} to $A$ of the representation $\tilde{\pi}$.
\end{defn}

\begin{prop}\label{decomposition}
Let the notation be as above and consider the restriction to $A$ of the
irreducible representations of $\tilde{A}$.  We then have:
$$
\tilde{\pi}_{r,s}|_A=\pi_r\qquad\text{and}\qquad\tilde{\pi}_{p,q}|_A=
\int^{\oplus}_{\mathbb{R}}\pi_{e^wp,e^{-w}q}\,dw.
$$
Here $\int^{\oplus}_{\mathbb{R}}$ denotes the direct integral (\cite{D1})
of representations.
\end{prop}

\begin{proof}
For any $\xi\in\tilde{\mathcal H}_{r,s}=L^2(\mathbb{R}^n)$, we have:
$$
\bigl(\tilde{Q}_{r,s}|_E(x,y;\theta)\xi\bigr)(u)=\theta\bar{e}\bigl[\eta_{\lambda}(r)
\beta(u,y)\bigr]\xi(u+x)=\bigl(Q_r(x,y;\theta)\xi\bigr)(u).
$$
It follows that $\tilde{\pi}_{r,s}|_A=\pi_r$.  Next, for any $\zeta\in\tilde
{\mathcal H}_{p,q}=L^2(\mathbb{R})$,
$$
\bigl(\tilde{Q}_{p,q}|_E(x,y;\theta)\zeta\bigr)(w)=\bar{e}(e^wp\cdot x+e^{-w}q\cdot y)
\zeta(w)=\bigl(Q_{e^wp,e^{-w}q}(x,y;\theta)\zeta\bigr)(w).
$$
By definition of the direct integrals, we thus obtain:
$$
\tilde{Q}_{p,q}|_E=\int^{\oplus}_{\mathbb{R}}Q_{e^wp,e^{-w}q}\,dw.
$$
It follows that: $\tilde{\pi}_{p,q}|_E=\int^{\oplus}_{\mathbb{R}}
\pi_{e^wp,e^{-w}q}\,dw$.
\end{proof}

\begin{rem}
This result has to do with the fact that $E$ is a normal subgroup of
$\tilde{E}$ with codimension 1.  Compare this result with Theorem 6.1 of
Kirillov's fundamental paper \cite{Ki2} or the discussion in section 2.5
of \cite{CG}, where the analysis of coadjoint orbits was used to obtain
a similar result for the representations of ordinary Lie groups.  Although
we proved this proposition directly, this strongly indicates the possibility
of formulating the proposition via generalized (dressing) orbit theory.
\end{rem}

 So far we have considered only the representations of $A$ and $\tilde{A}$.
In this article, we are not going to discuss the corepresentations of $(A,
\Delta)$ (or of $(\tilde{A},\tilde{\Delta})$).  In fact, it turns out that
the corepresentation theory is equivalent to the representation theory of
the dual group $G$.  This is a by-product of the Hopf $C^*$-algebra duality
between $(A,\Delta)$ and $(C^*(G),\hat{\Delta})$, provided by the ``regular''
multiplicative unitary operator $U$ associated with $(A,\Delta)$ (The
definition of $\hat{\Delta}$ depends on $U$.).  Since this is the case,
the corepresentation theory of $(A,\Delta)$ is actually simpler.

\section{Inner tensor product of representations}

 Given any two representations of a Hopf $C^*$-algebra $(B,\Delta)$, we 
can define their ``(inner) tensor product'' \cite[\S10]{Dr}, \cite[\S5]{CP}
as in the below.  There is also a corresponding notion for corepresentations.
But in the present article, we will not consider this dual notion.

\begin{defn}\label{tensor-product}
Let $\pi$ and $\rho$ be representations of a Hopf $C^*$-algebra $(B,\Delta)$,
acting on the Hilbert spaces ${\mathcal H}_{\pi}$ and ${\mathcal H}_{\rho}$.
Then their {\em inner tensor product\/}, denoted by $\pi\boxtimes\rho$, is a
representation of $B$ on ${\mathcal H}_{\pi}\otimes {\mathcal H}_{\rho}$
defined by
$$
(\pi\boxtimes\rho)(b)=(\pi\otimes\rho)\bigl(\Delta(b)\bigr),\qquad 
b\in B.
$$
Here $\pi\otimes\rho$ denotes the (outer) tensor product of the representations
$\pi$ and $\rho$, which is a representation of $B\otimes B$ naturally extended
to $M(B\otimes B)$.
\end{defn}

 As the name suggests, this notion of inner tensor product is a generalization
of the inner tensor product of group representations \cite{FD}.  For instance,
in the case of an ordinary group $C^*$-algebra $C^*(G)$ equipped with its
cocommutative (symmetric) comultiplication $\Delta_0$, Definition \ref
{tensor-product} is just the integrated form version of the inner tensor product
group representations.  In this case, since $\Delta_0$ is cocommutative, the
flip $\sigma:{\mathcal H}_{\pi}\otimes{\mathcal H}_{\rho}\to{\mathcal H}_{\rho}
\otimes{\mathcal H}_{\pi}$ provides a natural intertwining operator between
$\pi\boxtimes\rho$ and $\rho\boxtimes\pi$.

 For general (non-cocommutative) Hopf $C^*$-algebras, however, this is not
necessarily true.  In general, $\pi\boxtimes\rho$ need not even be equivalent
to $\rho\boxtimes\pi$.  Thus for any Hopf $C^*$-algebra (or quantum group),
it is an interesting question to ask whether two inner tensor product
representations $\pi\boxtimes\rho$ and $\rho\boxtimes\pi$ are equivalent
and if so, what the intertwining unitary operator between them is.

 When a Hopf $C^*$-algebra is equipped with a certain ``quantum universal
$R$-matrix'' (\cite{Dr}, \cite{CP}, and Definitions 6.1, 6.2 of \cite{BJKp2}
for the $C^*$-algebra version), we can give positive answers to these questions.
The following result is relatively well known:

\begin{prop}\label{equivitp}
Let $(B,\Delta)$ be a Hopf $C^*$-algebra.  Suppose that there exists a
quantum universal $R$-matrix $R\in M(B\otimes B)$ for $B$.  Then given any
two representations $\pi$ and $\rho$ of $B$ on Hilbert spaces ${\mathcal H}_{\pi}$
and ${\mathcal H}_{\rho}$, their inner tensor products $\pi\boxtimes\rho$ and
$\rho\boxtimes\pi$ are equivalent.  The equivalence is established by the
(unitary) intertwining operator $T_{\pi\rho}:{\mathcal H}_{\pi}\otimes
{\mathcal H}_{\rho}\to{\mathcal H}_{\rho}\otimes{\mathcal H}_{\pi}$, defined by
$$
T_{\pi\rho}=\sigma\circ(\pi\otimes\rho)(R).
$$
Here $\pi\otimes\rho$ is understood as the natural extension to $M(B\otimes B)$
of the tensor product $\pi\otimes\rho:B\otimes B\to{\mathcal B}({\mathcal H}_{\pi}
\otimes{\mathcal H}_{\rho})$ and $\sigma$ is the flip.  That is,
$$
T_{\pi\rho}\bigl((\pi\boxtimes\rho)(b)\bigr)=\bigl((\rho\boxtimes\pi)(b)\bigr)
T_{\pi\rho},\qquad b\in B.
$$
Furthermore, if the $R$-matrix is triangular, then we also have:
$$
T_{\rho\pi}T_{\pi\rho}=I_{{\mathcal H}_{\pi}\otimes{\mathcal H}_{\rho}}
\quad\text{and}\quad T_{\pi\rho}T_{\rho\pi}=I_{{\mathcal H}_{\rho}
\otimes{\mathcal H}_{\pi}}.
$$
\end{prop}

\begin{proof}
Let us first calculate how $T_{\pi\rho}$ acts as an operator.  If $\zeta
_1\in{\mathcal H}_{\pi}$ and $\zeta_2\in{\mathcal H}_{\rho}$, we have:
$$
T_{\pi\rho}(\zeta_1\otimes\zeta_2)=\sigma\circ\bigl((\pi\otimes\rho)(R)
\bigr)(\zeta_1\otimes\zeta_2)=\bigl((\rho\otimes\pi)(R_{21})\bigr)
(\zeta_2\otimes\zeta_1).
$$
Note that we have $T_{\pi\rho}=\bigl((\rho\otimes\pi)(R_{21})\bigr)
\circ\sigma$, as an operator.  To verify that $T_{\pi\rho}$ is an 
intertwining operator between $\pi\boxtimes\rho$ and $\rho\boxtimes\pi$,
let us consider an element $b\in B$.  Then
\begin{align}
T_{\pi\rho}\bigl((\pi\boxtimes\rho)(b)\bigr)&=T_{\pi\rho}\bigl
((\pi\otimes\rho)(\Delta b)\bigr)=\bigl((\rho\otimes\pi)(R_{21})\bigr)
\bigl((\rho\otimes\pi)(\Delta ^{\text{op}}b)\bigr)\circ\sigma  \notag \\
&=\bigl((\rho\otimes\pi)(R_{21}\Delta ^{\text {op}}b)\bigr)\circ\sigma.  \notag
\end{align}
From $R\Delta(b)R^{-1}=\Delta^{\text {op}}(b)$, we have: $R_{21}\Delta
^{\text {op}}(b)R_{21}^{-1}=\Delta(b)$.  It follows that
\begin{align}
T_{\pi\rho}\bigl((\pi\boxtimes\rho)(b)\bigr)&
=\bigl((\rho\otimes\pi)(\Delta(b)R_{21})\bigr)\circ\sigma
=\bigl((\rho\otimes\pi)(\Delta b)\bigr)\bigl((\rho\otimes\pi)(R_{21})
\bigr)\circ\sigma  \notag \\
&=\bigl((\rho\boxtimes\pi)(b)\bigr)T_{\pi\rho}.  \notag
\end{align}

Finally, if $R$ is triangular, by definition $\sigma\circ R=R_{21}=R^{-1}$.
Therefore,
\begin{align}
T_{\rho\pi}T_{\pi\rho}(\zeta_1\otimes\zeta_2)&=\bigl((\pi\otimes\rho)
(R_{21})\bigr)\bigl((\pi\otimes\rho)(R)\bigr)(\zeta_1\otimes\zeta_2) \notag \\
&=\bigl((\pi\otimes\rho)(R_{21}R)\bigr)(\zeta_1\otimes\zeta_2)
=(\zeta_1\otimes\zeta_2).  \notag
\end{align}
Since this is true for arbitrary $\zeta_1\in{\mathcal H}_{\pi}$ and 
$\zeta_2\in{\mathcal H}_{\rho}$, we have: $T_{\rho\pi}T_{\pi\rho}=I_
{{\mathcal H}_{\pi}\otimes{\mathcal H}_{\rho}}$.
\end{proof}

 In \cite[\S6]{BJKp2}, we showed that the ``extended'' Hopf $C^*$-algebra
$(\tilde{A},\tilde{\Delta})$ has a quasitriangular quantum universal $R$-matrix
$R\in M(\tilde{A}\otimes\tilde{A})$.  From this the following Corollary is immediate.

\begin{cor}
For the Hopf $C^*$-algebra $(\tilde{A},\tilde{\Delta})$, any two representations
$\tilde{\pi}$ and $\tilde{\rho}$ of $\tilde{A}$ will satisfy:
$$
\tilde{\pi}\boxtimes\tilde{\rho}\cong\tilde{\rho}\boxtimes\tilde{\pi}.
$$
By Proposition \ref{equivitp}, the operator $T_{\tilde{\pi}\tilde{\rho}}=
\sigma\circ(\tilde{\pi}\otimes\tilde{\rho})(R)$ is an intertwining operator
for this equivalence.
\end{cor}

 Unlike $(\tilde{A},\tilde{\Delta})$, however, the Hopf $C^*$-algebra $(A,\Delta)$
does not have its own quantum $R$-matrix $R_A\in M(A\otimes A)$.  Even at the
classical, Lie bialgebra level (studied in \cite[\S1]{BJKp2}), we can see that
the Poisson structures we consider cannot be obtained from any classical $r$-matrix.
Because of this, the result like the above Corollary is not automatic for
$(A,\Delta)$.  Even so, we plan to show in the below that the representations
of $(A,\Delta)$ still satisfy the quasitriangular type property.

 For this purpose and for possible future use, we are going to calculate here
the inner tensor product representations of our Hopf $C^*$-algebra $(A,\Delta)$.
Since it is sufficient to consider the inner tensor products of irreducible
representations, let us keep the notation of the previous section and let
$\{\pi_{p,q}\}_{(p,q)\in\mathbb{R}^{2n}}$ and $\{\pi_r\}_{r\in\mathbb{R}}$ be
the irreducible representations of $A$.  Similarly, let $\{\tilde{\pi}_s\}_
{s\in\mathbb{R}}$, $\{\tilde{\pi}_{p,q}\}_{(p,q)\in\mathbb{R}^{2n}}$,
$\{\pi_{r,s}\}_{(r,s)\in\mathbb{R}^2}$ be the irreducible representations
of $\tilde{A}$.  For convenience, we will calculate the inner tensor product
representations at the level of our dense subalgebra of functions, ${\mathcal A}$.

 Let $f\in{\mathcal A}$ and consider $\Delta f$.  To carry out our calculations,
it is convenient to regard $\Delta f$ also as a continuous function.  By using
the definition of $\Delta$ (given in Theorem 3.2 of \cite{BJKp2}) and by using
Fourier transform purely formally with the Fourier inversion theorem, it is not
difficult to realize $\Delta f$ as a function in the $(x,y,r)$ variables:
\begin{align}
&\Delta f(x,y,r,x',y',r')  \notag\\
&=\int f(x',y',r+r')\bar{e}\bigl[p\cdot(e^{\lambda r'}x'-x)+q\cdot
(e^{\lambda r'}y'-y)\bigr]\,dpdq. \notag
\end{align}
Note that in the $(p,q,r)\,(\in G)$ variables, it is just:
$$
\Delta f(p,q,r,p',q',r')=f(e^{\lambda r'}p+p',e^{\lambda r'}q+q',r+r'),
$$
which more or less reflects the multiplication law on $G$.

 We can now explicitly calculate the inner tensor products of irreducible
representations of ${\mathcal A}$, by $(\pi\boxtimes\rho)(f)=(\pi\otimes\rho)
(\Delta f)$.  We first begin with 1-dimensional representations.

\begin{prop}
For two 1-dimensional representations $\pi_{p,q}$ and $\pi_{p',q'}$
of $A$, we have: $\pi_{p,q}\boxtimes\pi_{p',q'}=\pi_{p+p',q+q'}$.  
From this, it follows that:
$$
\pi_{p,q}\boxtimes\pi_{p',q'}=\pi_{p+p',q+q'}=\pi_{p',q'}\boxtimes
\pi_{p,q}.
$$
\end{prop}

\begin{proof}
We have for any $f\in{\mathcal A}$,
\begin{align}
(\pi_{p,q}\boxtimes\pi_{p',q'})(f)&=\int f(x',y',0)\bar{e}\bigl[\tilde{p}
\cdot(x'-x)+\tilde{q}\cdot(y'-y)\bigr] \notag \\
&\qquad \bar{e}[p\cdot x+q\cdot y]\bar{e}[p'\cdot x'+q'\cdot y']
\,d\tilde{p}d\tilde{q}dxdydx'dy' \notag \\
&=\int f(x,y,0)\bar{e}\bigl[(p+p')\cdot x+(q+q')\cdot y\bigr]\,dxdy 
\notag \\ 
&=\pi_{p+p',q+q'}(f). \notag
\end{align}
\end{proof}

 For other cases involving infinite dimensional (irreducible) representations,
the equivalence between the inner tensor products is not so apparent.  However,
the inner tensor products of $\pi_{p,q}$ and $\pi_{r}$ has a property of being
equivalent to the infinite dimensional representation $\pi_r$ itself.  So in
this case, equivalence between the inner tensor products follows rather easily.

\begin{prop}
Consider a 1-dimensional representation $\pi_{p,q}$ and an infinite
dimensional representation $\pi_r$ of $A$.  Their inner tensor product
is equivalent to the irreducible representation $\pi_r$.  We thus have:
$$
\pi_r\boxtimes\pi_{p,q}\cong\pi_r\cong\pi_{p,q}\boxtimes\pi_r.
$$
\end{prop}

\begin{proof}
For $\xi\in{\mathcal H}_r\otimes{\mathcal H}_{p,q}\cong L^2(\mathbb{R}^n)
\otimes\mathbb{C}\cong L^2(\mathbb{R}^n)$ and for $f\in{\mathcal A}$, we
have:
\begin{align}
\bigl((\pi_r\boxtimes\pi_{p,q})(f)\xi\bigr)(u)&=\int f(x',y',r)\bar{e}
\bigl[\tilde{p}\cdot(x'-x)+\tilde{q}\cdot(y'-y)\bigr] \notag \\
&\qquad\bar{e}\bigl[\eta_{\lambda}(r)\beta(u,y)\bigr]\bar{e}[p\cdot x'
+q\cdot y']\xi(u+x)\,d\tilde{p}d\tilde{q}dxdydx'dy' \notag \\
&=\int f(x,y,r)\bar{e}[p\cdot x+q\cdot y]\bar{e}\bigr[\eta_{\lambda}(r)
\beta(u,y)\bigr]\xi(u+x)\,dxdy. \notag
\end{align}
Similarly for $\xi\in{\mathcal H}_{p,q}\otimes{\mathcal H}_r\cong
L^2(\mathbb{R}^n)$,
\begin{align}
\bigl((\pi_{p,q}\boxtimes\pi_r)(f)\xi\bigr)(u)&=\int f(x',y',r)\bar{e}
\bigl[\tilde{p}\cdot(e^{\lambda r}x'-x)+\tilde{q}\cdot(e^{\lambda r}y'-y)\bigr] 
\notag \\
&\qquad\bar{e}[p\cdot x+q\cdot y]\bar{e}\bigr[\eta_{\lambda}(r)
\beta(u,y')\bigr]\xi(u+x')\,d\tilde{p}d\tilde{q}dxdydx'dy' \notag \\
&=\int f(x,y,r)\bar{e}\bigl[e^{\lambda r}p\cdot x+e^{\lambda r}q\cdot y
\bigr]\bar{e}\bigl[\eta_{\lambda}(r)\beta(u,y)\bigr]\xi(u+x)\,dxdy. 
\notag
\end{align}

 Before proving $\pi_r\boxtimes\pi_{p,q}\cong\pi_{p,q}\boxtimes\pi_r$,
let us first show the equivalence $\pi_r\boxtimes\pi_{p,q}\cong\pi_r$.
This equivalence is suggested by the corresponding result at the level
of Heisenberg Lie group representation theory, which is obtained by
using the standard analysis via ``characters'' \cite{Ki}, \cite{CG}.
In our case, the equivalence is established by the intertwining operator
$S:L^2(\mathbb{R}^n)\to L^2(\mathbb{R}^n)$ defined by: $S\xi(u)=\bar{e}
(p\cdot u)\xi\left(u-\frac{q}{\eta_{\lambda}(r)}\right)$.  Indeed for
$f\in{\mathcal A}$,
\begin{align}
S\bigl((\pi_r\boxtimes\pi_{p,q})(f)\bigr)\xi(u)&=\int \bar{e}(p\cdot u)
f(x,y,r)\bar{e}[p\cdot x+q\cdot y]  \notag \\
&\qquad\bar{e}\left[\eta_{\lambda}(r)\beta\left(u-\frac{q}{\eta_{\lambda}(r)},
y\right)\right]\xi\left(u-\frac{q}{\eta_{\lambda}(r)}+x\right)\,dxdy, \notag
\end{align}
and
$$
\bigl(\pi_r(f)\bigr)S\xi(u)=\int f(x,y,r)\bar{e}\bigl[\eta_{\lambda}(r)
\beta(u,y)\bigr]\bar{e}\bigl[p\cdot(u+x)\bigr]\xi\left(u+x-\frac{q}
{\eta_{\lambda}(r)}\right)\,dxdy.
$$
We thus have: $S\bigl((\pi_r\boxtimes\pi_{p,q})(f)\bigr)=\bigl(\pi_r(f)
\bigr)S$, proving the equivalence: $\pi_r\boxtimes\pi_{p,q}\cong\pi_r$.
It is easy to check that $S^{-1}$ gives the intertwining operator for the
equivalence: $\pi_{r}\cong\pi_r\boxtimes\pi_{p,q}$.

 Meanwhile, from the explicit calculations given at the beginning of the
proof, it is apparent that we have: $\pi_{p,q}\boxtimes\pi_r=\pi_r\boxtimes
\pi_{e^{\lambda r}p,e^{\lambda r}q}$.  We thus obtain the equivalence:
$\pi_{p,q}\boxtimes\pi_r\cong\pi_r$, via the intertwining operator similar
to the above $S$, replacing $p$ and $q$ with $e^{\lambda r}p$ and
$e^{\lambda r}q$.  Combining these results, we can find the intertwining
operator $T:L^2(\mathbb{R}^n)\to L^2(\mathbb{R}^n)$ between $\pi_r
\boxtimes\pi_{p,q}$ and $\pi_{p,q}\boxtimes\pi_r$, obtained by multiplying
the respective intertwining operators for the equivalences $\pi_r\boxtimes
\pi_{p,q}\cong\pi_r$ and $\pi_r\cong\pi_{p,q}\boxtimes\pi_r$.  By
straightforward calculation, we have the following expression for $T$:
$$
T\xi(u)=\bar{e}\left[p\cdot u-e^{\lambda r}p\cdot\left(u-\frac{q}
{\eta_{\lambda}(r)}+\frac{e^{\lambda r}q}{\eta_{\lambda}(r)}\right)
\right]\xi\left(u-\frac{q}{\eta_{\lambda}(r)}+\frac{e^{\lambda r}q}
{\eta_{\lambda}(r)}\right).
$$
It is clear that $T^{-1}$ gives the intertwining operator between $\pi_{p,q}
\boxtimes\pi_r$ and $\pi_r\boxtimes\pi_{p,q}$.
\end{proof}

 So far nothing very interesting has happened, in the sense that the results
are similar to those of the Heisenberg group representation theory.  However,
a breakdown of this analogy occurs when we consider inner tensor products of
two infinite dimensional representations $\pi_r$ and $\pi_{r'}$.  Let us first
prove the equivalence of $\pi_r\boxtimes\pi_{r'}$ and $\pi_{r'}\boxtimes\pi_{r}$.

\begin{prop}
Consider a pair $(\pi_r,\pi_r')$ of two infinite dimensional irreducible 
representations of $A$.  Then we have:
$$
\pi_r\boxtimes\pi_{r'}\cong\pi_{r'}\boxtimes\pi_r,
$$
where the equivalence between them is given by the intertwining operator
$T_{\pi_r\pi_{r'}}:L^2(\mathbb{R}^{2n})\to L^2(\mathbb{R}^{2n})$ defined by
$$
T_{\pi_r\pi_{r'}}\xi(v,u)=(e^{\frac{-\lambda r}2})^n(e^{-\frac
{\lambda r'}2})^n\,\xi\bigl(e^{-\lambda r'}u+(e^{\lambda r'}
-e^{-\lambda r'})e^{-\lambda r}v,e^{-\lambda r}v\bigr).
$$
\end{prop}

\begin{proof}
Let $f\in{\mathcal A}$.  For $\xi\in{\mathcal H}_r\otimes{\mathcal H}_{r'}
\cong L^2(\mathbb{R}^{2n})$, we have:
\begin{align}
&\bigl((\pi_r\boxtimes\pi_{r'})(f)\xi\bigr)(u,v)  \notag \\
&=\int f(x',y',r+r')\bar{e}\bigl[\tilde{p}\cdot(e^{\lambda r'}x'-x)
+\tilde{q}\cdot(e^{\lambda r'}y'-y)\bigr]  \notag \\
&\qquad\bar{e}\bigr[\eta_{\lambda}(r)\beta(u,y)\bigr]\bar{e}\bigr[
\eta_{\lambda}(r')\beta(v,y')\bigr]\xi(u+x,v+x')
\,d\tilde{p}d\tilde{q}dxdydx'dy' \notag \\
&=\int f(x,y,r+r')\bar{e}\bigl[\eta_{\lambda}(r)\beta(u,e^{\lambda r'}y)
\bigr]\bar{e}\bigl[\eta_{\lambda}(r')\beta(v,y)\bigr]
\xi(u+e^{\lambda r'}x,v+x)\,dxdy. \notag
\end{align}
Similarly for $\xi\in{\mathcal H}_{r'}\otimes{\mathcal H}_{r}\cong 
L^2(\mathbb{R}^{2n})$, by interchanging the roles of $r$ and $r'$,
\begin{align}
\bigl((\pi_{r'}\boxtimes\pi_{r})(f)\xi\bigr)(v,u)&=\int f(x,y,r+r')\bar{e}
\bigl[\eta_{\lambda}(r')\beta(v,e^{\lambda r}y)\bigr]\bar{e}\bigl[\eta_{\lambda}(r)
\beta(u,y)\bigr]  \notag \\
&\qquad\xi(v+e^{\lambda r}x,u+x)\,dxdy. \notag
\end{align}

 To prove the equivalence between $\pi_{r}\boxtimes\pi_{r'}$ and $\pi_{r'}
\boxtimes\pi_{r}$, it is useful to recall the fact that any two representations
$\tilde{\pi}$ and $\tilde{\rho}$ of the ``extended'' Hopf $C^*$--algebra
$(\tilde{A},\tilde{\Delta})$ satisfy: $\tilde{\pi}\boxtimes\tilde{\rho}\cong
\tilde{\rho}\boxtimes\tilde{\pi}$.  In particular, we would have:
$\tilde{\pi}_{r,0}\boxtimes\tilde{\pi}_{r',0}\cong\tilde{\pi}_{r',0}
\boxtimes\tilde{\pi}_{r,0}$.  Its intertwining operator is:
$T_{\tilde{\pi}_{r,0}\tilde{\pi}_{r',0}}=\sigma\circ\bigl((\tilde{\pi}_{r,0}
\otimes\tilde{\pi}_{r',0})(R)\bigr)$, by Corollary of Proposition \ref{equivitp}.
By restriction to $A$ (in the sense of Definition \ref{restriction}), we obtain:
$$
(\tilde{\pi}_{r,0}|_A)\boxtimes(\tilde{\pi}_{r',0}|_A)\cong(\tilde{\pi}
_{r',0}|_A)\boxtimes(\tilde{\pi}_{r,0}|_A),
$$
which, by Proposition \ref{decomposition}, is just:
$\pi_{r}\boxtimes\pi_{r'}\cong\pi_{r'}\boxtimes\pi_{r}$.

 To find the intertwining operator for this equivalence, let us find an
explicit expression for the operator $T_{\tilde{\pi}_{r,0}\tilde{\pi}_{r',0}}$.
Recall first that by equation (6.3) and Definition 6.3 of \cite{BJKp2}, the
quantum $R$-matrix for $(\tilde{A},\tilde{\Delta})$ is considered as a continuous
function defined by:
\begin{align}
R(p,q,r,s,p',q',r',s')&=\Phi(p,q,r,s,p',q',r',s')\Phi'(p,q,r,s,p',q',r',s')
\notag \\
&=\bar{e}\bigl[\lambda(rs'+r's)\bigr]\bar{e}\bigl[2\lambda(e^{-\lambda r'})
p\cdot q'\bigr].  \notag
\end{align}
We then calculate $(\tilde{\pi}_{r,0}\otimes\tilde{\pi}_{r',0})(R)$ as an
operator on $\tilde{\mathcal H}_{r,0}\otimes\tilde{\mathcal H}_{r',0}\cong
L^2(\mathbb{R}^{2n})$.  By a straightforward calculation, we obtain for
$\xi\in L^2(\mathbb{R}^{2n})$,
\begin{align}
&\bigl((\tilde{\pi}_{r,0}\otimes\tilde{\pi}_{r',0})(\Phi)\xi\bigr)(u,v)
\notag \\
&=\int \bar{e}\bigl[\lambda(rs'+r's)\bigr]e[p\cdot x+q\cdot y+sw+p'\cdot
x'+q'\cdot y'+s'w']  \notag \\
&\qquad(e^{-\frac{w}2})^n(e^{-\frac{w'}2})^n\,\xi(e^{-w}u+e^{-w}x,
e^{-w'}v+e^{-w'}x') \notag \\
&\qquad\bar{e}\bigl[\eta_{\lambda}(r)\beta(u,y)\bigr]\bar{e}
\bigl[\eta_{\lambda}(r')\beta(v,y')\bigr]\,dpdqdsdp'dq'ds'dxdydwdx'dy'dw'
\notag \\
&=(e^{-\frac{\lambda r}2})^n(e^{-\frac{\lambda r'}2})^n\,\xi(e^{-\lambda r'}
u,e^{-\lambda r}v). \notag
\end{align}
Similarly,
$$
\bigl((\tilde{\pi}_{r,0}\otimes\tilde{\pi}_{r',0})(\Psi)\xi\bigr)(u,v)
=\xi(u+2\lambda e^{-\lambda r'}\eta_{\lambda}(r')v,v).
$$
Since $\eta_{\lambda}(r')=\frac{e^{2\lambda r'}-1}{2\lambda}$, we thus have:
\begin{align}
\bigl((\tilde{\pi}_{r,0}\otimes\tilde{\pi}_{r',0})(R)\xi\bigr)(u,v)
&=\bigl((\tilde{\pi}_{r,0}\otimes\tilde{\pi}_{r',0})(\Phi)\bigr)
\bigl((\tilde{\pi}_{r,0}\otimes\tilde{\pi}_{r',0})(\Psi)\bigr)\xi(u,v)
\notag \\
&=(e^{\frac{-\lambda r}2})^n(e^{-\frac{\lambda r'}2})^n\xi
\bigl(e^{-\lambda r'}u+(e^{\lambda r'}-e^{-\lambda r'})e^{-\lambda r}v,
e^{-\lambda r}v\bigr). \notag
\end{align}

 By applying the flip $\sigma$, we therefore obtain:
\begin{align}
T_{\tilde{\pi}_{r,0}\tilde{\pi}_{r',0}}\xi(v,u)&=\sigma\circ\bigl((\tilde
{\pi}_{r,0}\otimes\tilde{\pi}_{r',0})(R)\bigr)\xi(v,u)  \notag \\
&=(e^{\frac{-\lambda r}2})^n(e^{-\frac{\lambda r'}2})^n\xi
\bigl(e^{-\lambda r'}u+(e^{\lambda r'}-e^{-\lambda r'})e^{-\lambda r}v,
e^{-\lambda r}v\bigr). \notag
\end{align}

 Define $T_{\pi_r\pi_{r'}}$ by $T_{\pi_r\pi_{r'}}=T_{\tilde{\pi}_{r,0}
\tilde{\pi}_{r',0}}$.  Then it is a straightforward calculation to show
that $T_{\pi_r\pi_{r'}}$ is an intertwining operator between $\pi_{r}\boxtimes
\pi_{r'}$ and $\pi_{r'}\boxtimes\pi_{r}$.  For $f\in{\mathcal A}$ and
$\xi\in L^2(\mathbb{R}^{2n})$, we have:
\small
\begin{align}
&T_{\pi_r\pi_{r'}}\bigl((\pi_{r}\boxtimes\pi_{r'})(f)\bigr)\xi(v,u) 
\notag \\
&=\int f(x,y,r+r')\bar{e}\bigl[\eta_{\lambda}(r)\beta(e^{-\lambda r'}u
+(e^{\lambda r'}-e^{-\lambda r'})e^{-\lambda r}v,e^{\lambda r'}y)\bigr]
\bar{e}\bigl[\eta_{\lambda}(r')\beta(e^{-\lambda r}v,y)\bigr] \notag \\
&\quad\quad\quad(e^{\frac{-\lambda r}2})^n(e^{-\frac{\lambda r'}2})^n
\xi(e^{-\lambda r'}u+(e^{\lambda r'}-e^{-\lambda r'})e^{-\lambda r}v
+e^{\lambda r'}x,e^{-\lambda r}v+x)\,dxdy \notag \\
&=\bigl((\pi_{r'}\boxtimes\pi_{r})(f)\bigr)T_{\pi_r\pi_{r'}}\xi(v,u).
\notag
\end{align}
\normalsize
Since $\xi$ is arbitrary, it follows that: 
$$
T_{\pi_r\pi_{r'}}\bigl((\pi_{r}\boxtimes\pi_{r'})(f)\bigr)=\bigl((\pi_{r'}
\boxtimes\pi_{r})(f)\bigr)T_{\pi_r\pi_{r'}}.
$$
\end{proof}

 Observe that by interchanging $r$ and $r'$, we are able to find the expression
for the intertwining operator $T_{\pi_{r'}\pi_{r}}$: 
$$
T_{\pi_{r'}\pi_{r}}\xi(u,v)=(e^{\frac{-\lambda r}2})^n(e^{-\frac
{\lambda r'}2})^n\,\xi\bigl(e^{-\lambda r}v+(e^{\lambda r}-e^{-\lambda r})
e^{-\lambda r'}u,e^{-\lambda r'}u\bigr).
$$
For $f\in{\mathcal A}$, we will have: $T_{\pi_r'\pi_{r}}\bigl((\pi_{r'}
\boxtimes\pi_{r})(f)\bigr)=\bigl((\pi_{r}\boxtimes\pi_{r'})(f)\bigr)
T_{\pi_r'\pi_{r}}$.

 However, note that unlike in the ordinary group representation theory or in
the cases equipped with ``triangular'' quantum $R$-matrices, we no longer have: 
$T_{\pi_{r'}\pi_{r}}T_{\pi_r\pi_{r'}}=I$.  We instead have:
\small
\begin{align}
T_{\pi_{r'}\pi_{r}}T_{\pi_{r}\pi_{r'}}\xi(u,v)=(e^{-\lambda r})^n
(e^{\lambda r'})^n\,\xi\bigl(&u-(1-e^{-2\lambda r'})e^{-2\lambda r}u+
(e^{\lambda r'}-e^{-\lambda r'})e^{-2\lambda r}v,  \notag \\
&\quad e^{-2\lambda r}v+
(1-e^{-2\lambda r})e^{-\lambda r'}u\bigr),  \notag
\end{align}
\normalsize
which is clearly not the identity operator.  Let us summarize our results
in the following theorem:
\begin{theorem}
Given any two representations $\pi$ and $\rho$ (acting on the Hilbert spaces
${\mathcal H}_{\pi}$ and ${\mathcal H}_{\rho}$) of our Hopf $C^*$-algebra
$(A,\Delta)$, their inner tensor products ``commutes'' (i.\,e. $\pi\boxtimes
\rho$ and $\rho\boxtimes\pi$ are equivalent).  However, the intertwining
operators between them behave in an interesting way, in the sense that we have:
$T_{\rho\pi}T_{\pi\rho}\ne I_{{\mathcal H}_{\pi}\otimes{\mathcal H}_{\rho}}$,
in general.
\end{theorem}

 This theorem means that the category of representations of $(A,\Delta)$
is essentially a ``quasitriangular monoidal category''.  This is a typical
characteristic of the category of ``braids'' (in knot theory).  In recent
years, representation theory of quantum groups led to the developments and
discoveries of some useful knot invariants, like Jones polynomials or HOMFLY
polynomials, especially in connection with the existence of quantum universal
$R$-matrices (For more discussion about these topics, see \cite{Bir} or
\cite[\S15]{CP}.).

 In our case, it is interesting to point out that $(A,\Delta)$ possesses
the quasitriangular type property, without the existence of its own quantum
$R$-matrix $R_A$.  Meanwhile, since there have been only a handful of examples
so far of non-compact, $C^*$-algebraic quantum groups possessing the property
of quasitriangularity, having these examples $(A,\Delta)$ and $(\tilde{A},
\tilde{\Delta})$ would benefit the study of non-compact quantum groups and
its development.



\bibliography{ref}

\bibliographystyle{amsplain}

\end{document}